\documentclass[12pt,leqno]{article}
\usepackage{amssymb,amsmath,amsthm}
\usepackage[all]{xy}
\numberwithin{equation}{section}
\date{}

\usepackage{amssymb}
\usepackage{amsbsy}
\usepackage{amscd}
\usepackage{amsmath}
\usepackage{amsthm}
\usepackage[mathscr]{eucal}

\usepackage{kyoto}
\usepackage[all]{xy}
\usepackage{graphicx}

\title{Masaki Kashiwara and Algebraic Analysis}

\author{Pierre Schapira}

\begin{document}

\maketitle
\begin{abstract}
This paper is based on a talk given in honor of Masaki Kashiwara's  
60th birthday, Kyoto, June 27, 2007.
It is a brief overview of his main contributions in the domain of 
microlocal and algebraic analysis.
\end{abstract}

It is a great honor
to present here some aspects of the 
work of Masaki Kashiwara.

Recall that Masaki's work covers many fields of
mathematics, 
algebraic and microlocal analysis of course,  
but also representation theory, 
Hodge theory, integrable systems, quantum groups and so on. 
Also recall that Masaki had many collaborators, among whom 
Daniel Barlet, Jean-Luc Brylinski, Etsurio Date, Ryoji Hotta, Michio Jimbo, 
Seok-Jin Kang, Takahiro Kawai,
Tetsuji Miwa, Kiyosato Okamoto, Toshio Oshima, Mikio Sato, myself, 
 Toshiyuki Tanisaki and Mich{\`e}le Vergne.

In each of the domain he approached, Masaki has
given essential contributions and made important discoveries, such
as, for example, the existence of crystal bases in quantum groups. 
But in this talk, I will restrict myself to describe some part of his
work related to microlocal and algebraic analysis.

The story begins long ago, in the early sixties, 
when Mikio Sato created a new branch of mathematics
 now called  ``Algebraic Analysis''.  In 1959/60, M.~Sato published
 two papers on hyperfunction theory \cite{Sa1} and then 
developed his vision of analysis and linear partial differential
equations in a series of
lectures at Tokyo University (see \cite{An}).
If $M$ is a real analytic manifold and $X$ is a complexification of
$M$, hyperfunctions on $M$ are cohomology classes supported by $M$ 
of the sheaf $\sho_X$ of holomorphic functions on $X$. 
It is difficult to realize now how Sato's point of view
was revolutionary at that time. 
Sato's hyperfunctions are constructed using tools from sheaf
theory and complex analysis, when people were totally
addicted to functional analysis, 
and when the separation between real and complex analysis was very strong.

\begin{figure}
\begin{tabular}{cc}
\includegraphics[scale=.3]{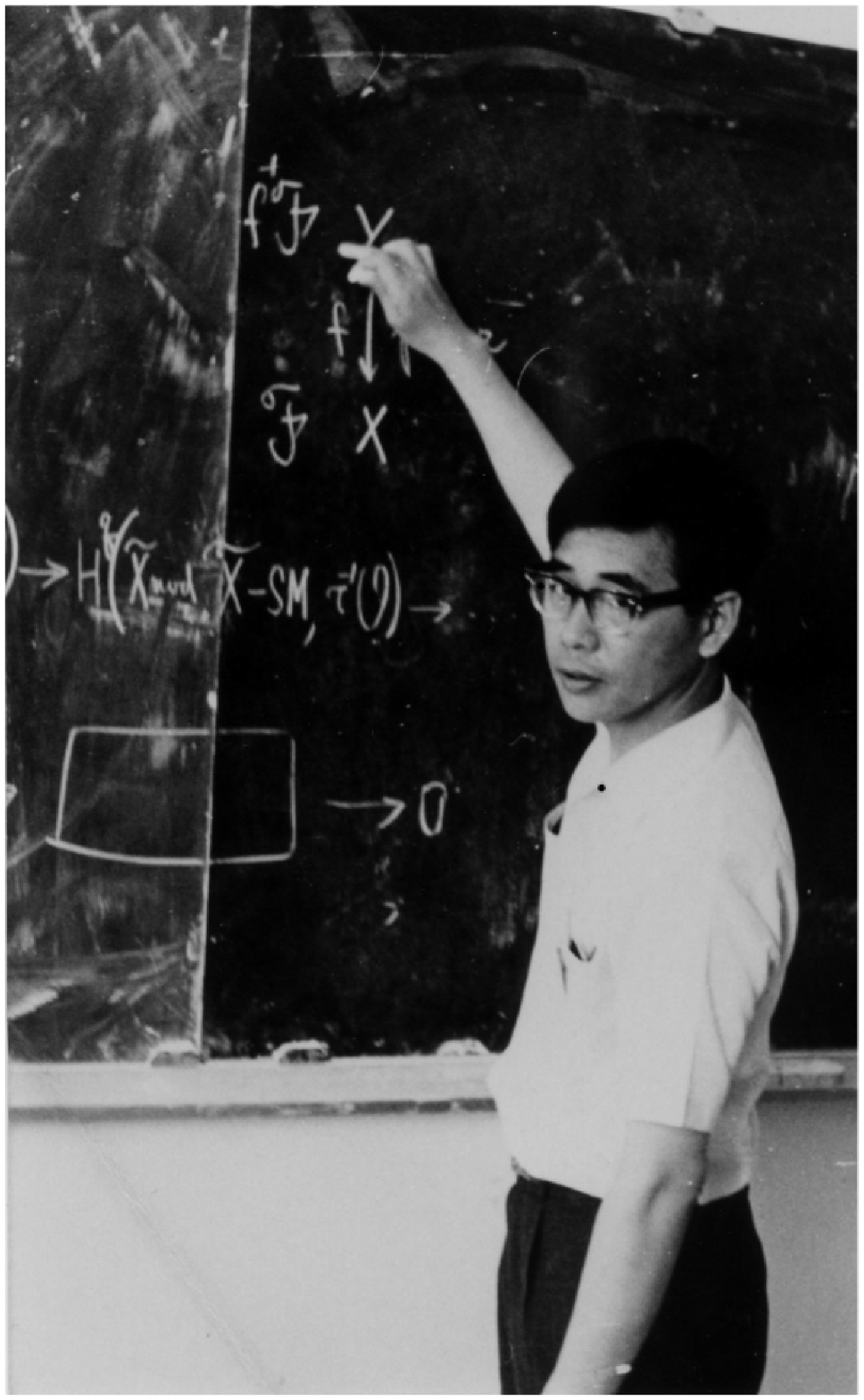}
&
\includegraphics[scale=.5]{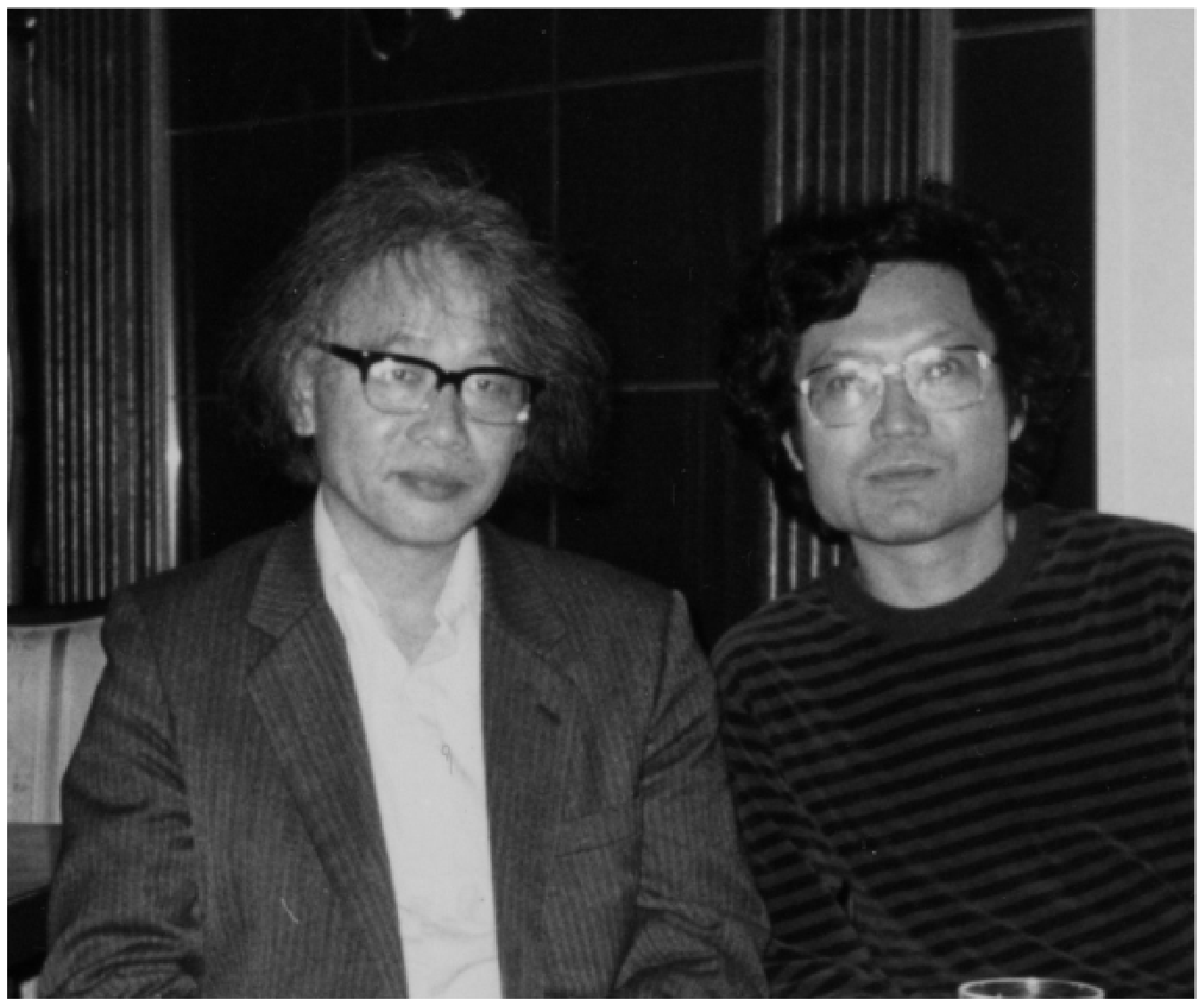}
\end{tabular}
\caption{Mikio Sato around 1972 and with Masaki Kashiwara more recently}
\end{figure}

Then came Kashiwara's thesis, dated December 1970 (of course
written in Japanese, but translated in English and published by the
French Mathematical Society \cite{Ka1}) in which he settles the foundations of analytic 
$\shd$-module theory and obtains almost all basic results of the theory
(compare with \cite{Ka7}). 
With $\shd$-module theory (also constructed independently in the algebraic
settings by J.~Bernstein \cite{Be}), 
one finally has the tools to treat general systems of
linear partial differential equations, as opposed to one equation with
one unknown, or to some very particular overdetermined systems. 
In particular, 
Kashiwara succeeds in formulating (and solving, but the difficult
problem is the functorial formulation) the Cauchy problem for $\shd$-modules, 
obtaining what is now called the Cauchy-Kowalesky-Kashiwara
theorem.

After the first revolution of hyperfunction theory, 
Sato made a second one, ten years later, by creating microlocal
analysis, a way to analyse objects of a manifold $X$ in the cotangent
bundle $T^*X$. 
With Kashiwara  and Kawai, they wrote a long paper \cite{S-K-K}, 
quoted everywhere as SKK, 
whose influence has been considerable during the whole seventies among 
the analysts (and not only the analysts), 
although very few of them even tried to read the paper. 
The SKK paper contains Sato's construction of the
sheaf $\shc_M$ of microfunctions, and as a byproduct, the definition of
the wave front set. This is essentially what the analysts,
led by H{\"o}rmander,  remember of this theory (see \cite{Ho}). 
But, to my opinion,  this is certainly not the only key point of
the SKK paper. Another essential fact is that all
constructions are made functorially. For example, microfunctions are
obtained by first constructing the microlocalization functor $\mu_M$, and
then applying it to the sheaf $\sho_X$ of holomorphic functions on a
complex manifold $X$. When you take for $M$ a real analytic manifold 
of whom $X$ is a complexification, you get the sheaf $\shc_M$ (living
on $T^*_MX$, the cornormal bundle to $M$ in $X$), but if 
you replace the embedding $M\hookrightarrow X$ by the diagonal
embedding $\Delta\hookrightarrow X\times X$, then you get
the sheaf of
microdifferential operators (on $T^*_\Delta (X\times X)\simeq T^*X$)
whose theory was developed by Kashiwara
and Kawai. With this approach, you can adapt the six Grothendieck
operations to Analysis 
and obtain a completely new point of view to classical
problems ({\em e.g.} the Fourier-Sato transformation). 

Moreover, the SKK paper contains at least two fundamental
and extremely deep results, 
first the involutivity of characteristics, second the 
structure of systems of microdifferential equations at the generic
points of the characteristic variety. More precisely, 
let  $X$ be a complex manifold and let $\she_X$ be the sheaf of rings of
microdifferential operators (a kind of localization of the sheaf $\shd_X$ of differential operators). A 
microdifferential system $\shm$ on an open subset 
$U$ of $T^*X$ is a coherent $\she_X\vert_U$-module. Then

\begin{itemize}
\item
the support $\chv(\shm)$ of $\shm$, also called its characteristic
variety, is a closed complex analytic {\em involutive} (that is, co-isotropic)
subset of $U$. Of course, the involutivity theorem has a longer history,
including the previous work of Guillemin-Quillen-Ster\-nberg \cite{GQS}, 
and culminating
with the purely algebraic proof of Gabber \cite{Ga}. 

\item
At {\em generic points} of $\chv(\shm)$, 
(after using complex quantized contact transformations and 
infinite order microdifferential operators) 
$\shm$ is isomorphic to a partial de Rham system:
\eqn
&&\partial_{x_i}u=0,\,(i=1,\dots,p).
\eneqn
In the real case, $\shm$ is isomorphic to a mixture of 
de Rham, Dolbeault and Hans Lewy systems:
\eqn
\left\{\begin{array}{l}
\partial_{x_i}u=0,\,(i=1,\dots,p)\\
(\partial_{y_j}+\sqrt{-1}\partial_{y_{j+1}})u=0,(j=1,\dots,q)\\
(\partial_{t_k}+\sqrt{-1}t_k\partial_{t_{k+1}})u=0,(k=1,\dots,r).
\end{array}\right.
\eneqn
\end{itemize}

From 1970 to 1980, Kashiwara solved almost all fundamental questions of 
$\shd$-module theory, proving in particular the rationality of the zeroes
of $b$-functions \cite{Ka4} and also stating and solving almost 
all questions related to
regular holonomic modules, in particular the Riemann-Hilbert problem.
\begin{figure}
\centerline{
\includegraphics[scale=.5]{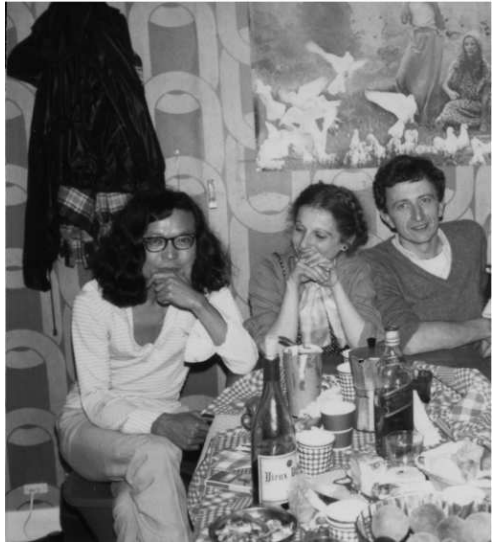}}
\caption{Masaki Kashiwara, Teresa Monteiro Fernandes and myself around 1975}
\end{figure}

Let us give some details on this part of Kashiwara's work. 
In 1975, he proved that
the complex
$F=\rhom[\shd](\shm,\sho_X)$ of holomorphic solutions of a holonomic 
$\shd$-module $\shm$ has constructible cohomology and satisfies 
properties which are now translated by 
saying that $F$ is perverse \cite{Ka3}. Moreover, 
two years {\it before} \cite{Ka2}, in 1973, he 
calculated the local 
Euler-Poincar{\'e} index of $F$ using the characteristic cycle
associated to $\shm$ and in fact, defining first what is now called 
the local Euler obstruction, or equivalently, 
the intersection of Lagrangian cycles.
In 1977 he gave a precise statement of what should be the Riemann-Hilbert 
correspondence (see \cite[p.~287]{Ra}), the difficulty  
being to define a suitable class of 
holonomic $\shd$-modules, the so-called regular holonomic modules, 
what he does in the microlocal
setting with Kawai \cite{Kaw} (after related work with Oshima \cite{Kao}).
Then, in 1979, he announces at the 
1979/1980 Seminar of Ecole Polytechnique \cite{Ka5} 
the theorem, giving with some details the main steps of the proof. 
\eqn
&&\xymatrix{\RD^\Rb_\hol(\shd_X)^{\rop}\ar[rr]^-{(1975)}
             &&\RD^\Rb_{\Cc}(\C_X)\ar@<0.5ex>[lld]_-\sim^-{(1979-80)}\\
\RD^\Rb_{\holreg}(\shd_X)^{\rop}\ar@{^(->}[u]^-{(1977)}&&
}
\eneqn
Unfortunately, Masaki did not publish the 
whole proof before 1984 \cite{Ka6} and some people tried to make 
his result their own. As everyone knows, 
if the platonic world of Mathematics is pure and rigorous, these
qualities definitely do not apply to the world of mathematicians.
 
Of course, Kashiwara did a lot of other things during this period 
1970/80, in particular in the theory of microdifferential equations,
but he did not always take the time to publish his results. 
I remember that I had once in 1978 at Oberwolfach  
the opportunity to explain to 
H{\"o}rmander the so-called ``watermelon cut theorem'' and you can now find
it in \cite[Th. 9.6.6]{Ho}. 
This beautiful theorem
asserts in particular that if a hyperfunction $u$  is supported by a half space
$f\geq 0$, then the
analytic wave front set of $u$ above the boundary $f=0$ is invariant by
the Hamiltonian vector field $H_f$. 

After that, essentially from 1980 to 1990, 
 came another period in which I am more involved. 

\begin{figure}
\centerline{
\includegraphics[scale=.4]{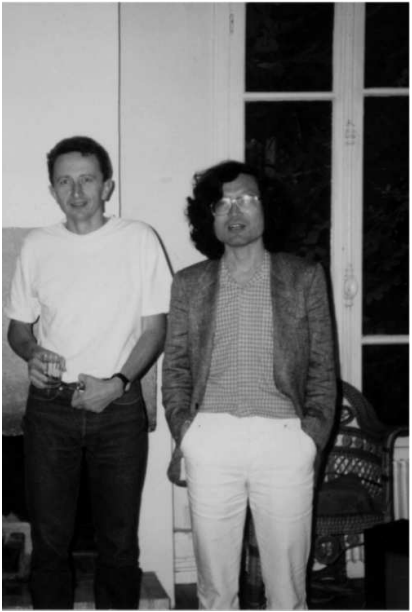}}
\caption{Masaki Kashiwara and myself around 1980}
\end{figure}
Indeed, we developed 
together the microlocal theory of sheaves (see \cite{KS1}). To a sheaf $F$
(not necessarily constructible) on a {\em real} manifold $M$, 
we associate a closed 
conic subset  $\SSu(F)$\footnote{$\SSu(F)$ stands for singular support.} 
of the cotangent bundle, the microsupport of $F$, 
which describes the directions of non propagation of $F$. 
The idea of microsupport emerged when, on one side,
Masaki noticed that it was possible to recover the
characteristic variety of a holonomic $\shd$-module from the knowledge
of the complex of its holomorphic solutions by using the vanishing
cycle functor, and when, on my side, I was lead to this notion by 
remarking that our previous results on  propagation for hyperbolic
systems  was of purely geometrical nature
and had almost nothing to do with partial differential equations.

One of the main result of the theory asserts that the 
microsupport of a sheaf $F$ is an involutive subset of the 
cotangent bundle, but now we are working on real manifolds.
In case one works on a complex manifold 
and $F$ is the sheaf of solutions of a coherent $\shd$-module $\shm$, 
the microsupport of $F$ coincides with the characteristic variety of $\shm$. 
This gives an alternative proof of the involutivity of
characteristics of $\shd$-modules. 
Moreover, constructible sheaves on a real manifold are sheaves whose
microsupport is subanalytic and Lagrangian.
This allowed  Kashiwara to 
adapt to the real case the notion of characteristic cycle of a
$\shd$-module and to 
define the Lagrangian cycle of an $\R$-constructible sheaf. 
The group of Lagrangian cycles 
is isomorphic to the Grothendieck group of the abelian category of 
$\R$-constructible sheaves and is also isomorphic to the group of 
constructible 
functions. Lagrangian cycles play a basic role in many questions and
have been recently extended to higher $K$-theory by Beilinson \cite{Bei}.  

After 90, Kashiwara concentrated mainly on other subjects such as
crystal bases, but nevertheless we wrote several papers together.
In order to overcome some difficulties
related to the microlocalization functor, 
we were led to generalize the notion
of sheaves and to define ind-sheaves \cite{KS2}. This theory required a lot 
of technology from category theory, and, as a byproduct, we wrote
a whole book on this subject \cite{KS3}. 

Algebraic Analysis and Microlocal Analysis are still actively developing in
various directions. Let us mention three of them. 

(i) Recall that Masaki was the first, in 1996, to introduce algebroid
stacks in microlocal analysis \cite{Ka8}. 
Indeed, on a complex contact manifold the sheaf of
microdifferential operators does not exist in general and one
has to replace sheaves with stacks. 
Such algebroid stacks are now commonly used on 
complex symplectic manifolds where microdifferential operators are
replaced by a variant involving a central parameter $\hbar$.
Note that Masaki and Raphael Rouquier 
 recently used such rings of operators to make  a surprising link 
with Cherednik algebras \cite{KaR}.

(ii) As a particular case of the theory of ind-sheaves,
 one gets the theory of usual sheaves on the subanalytic site. 
Personally, (I am not sure that Masaki shares this point of
view) I am convinced that the subanalytic topology is 
particularly well suited
to treat many problems in Analysis and that there are 
lot of interesting results to be obtained in this direction.

(iii) Another very promising direction is the link between the microlocal
theory of sheaves and Fukaya's category.
On one side,  Nadler and Zaslow \cite{NaZ,Na}, adapting the construction of 
Lagrangian cycles, constructed a category equivalent to Fukaya's category 
(on cotangent bundles, not  on general real symplectic manifolds).
On the other side,   
Tamarkin \cite{Ta} also constructed a category which should play
this role. 
Tamarkin's idea is to add a variable $t\in\R$ whose dual variable
$\tau$ plays the role of the inverse of $\hbar$ and to
 work ``microlocally'' with the  category of constructible sheaves 
on $X\times\R$ in the open set $\tau>0$ of $T^*(X\times\R)$. 

I hope that this very sketchy panorama of almost fifty years of 
Algebraic Analysis (perhaps one should now better call it 
``Functorial Analysis'')  
will have convinced you of the importance of the
theory and of the fact that Masaki plays the main role in it since 
the early seventies.

\providecommand{\bysame}{\leavevmode\hbox to3em{\hrulefill}\thinspace}

\vspace*{1cm}
\noindent
\parbox[t]{16em}
{\scriptsize{
Pierre Schapira\\
Institut de Math{\'e}matiques\\
Universit{\'e} Pierre et Marie Curie\\
175, rue du Chevaleret,
75013 Paris, France\\
e-mail: schapira@math.jussieu.fr\\
http://www.math.jussieu.fr/$\sim$schapira/
}}

\end{document}